\documentclass[letterpaper, 10pt, conference]{ieeeconf}

\IEEEoverridecommandlockouts
\usepackage{amssymb}
\usepackage{amsmath}
\usepackage{amsfonts}
\usepackage{graphicx}
\usepackage{epsfig}
\usepackage{textcomp}
\usepackage{upgreek}
\usepackage{bm}
\usepackage{xcolor}

\newtheorem{theorem}{Theorem}

\newtheorem{lemma}{Lemma}

\DeclareSymbolFont{symbolsC}{U}{pxsyc}{m}{n}
\DeclareMathSymbol{\originalimage}{\mathrel}{symbolsC}{23}
\DeclareMathSymbol{\imageoriginal}{\mathrel}{symbolsC}{24}
\DeclareMathSymbol{\imageoriginalV}{\mathrel}{symbolsC}{151}
\DeclareMathSymbol{\originalimageV}{\mathrel}{symbolsC}{152}

\begin{document}

\title{\bf \Large On algebraic time-derivative estimation \\
  and deadbeat state reconstruction}
\author{Johann Reger$^{\dag}$\thanks{$^{\dag }$Johann Reger is head of
    the Control Engineering Group, Ilmenau University of Technology,
    Gustav-Kirchhoff-Str.1, D-98693 Ilmenau, Germany (e-mail:
    reger@ieee.org).} and Jerome Jouffroy$^{\ddag}$\thanks{%
    $^{\ddag }$Jerome Jouffroy is with the Mads Clausen Institute,
    University of Southern Denmark, Alsion 2, DK-6400 S\o nderborg, Denmark
    (e-mail: jerome@mci.sdu.dk).}}

\bibliographystyle{IEEEtran}
\thispagestyle{empty}
\pagestyle{empty}

\maketitle%

\begin{abstract}%
This paper places into perspective the so-called algebraic time-derivative
estimation method recently introduced by Fliess and co-authors with standard
results from linear state-space theory for control systems. In particular,
it is shown that the algebraic method can essentially be seen as a special
case of deadbeat state estimation based on the reconstructibility Gramian of
the considered system.

\end{abstract}%

\section{Introduction}

In the past few years, the algebraic approach to estimation in control
systems proposed by Fliess and co-workers has generated a number of
interesting results for different problems of estimation of dynamical
systems such as state estimation, parametric identification, and fault
diagnosis, to name but a few (see \cite{Fliess-2004NOLCOS,
Fliess-2003esaim,Fliess-2007easy,Fliess-2004diagijc} and references
therein). Loosely speaking, this new estimation approach is mainly based on
the robust computation of the time-derivatives of a noisy signal by using a
finite weighted combination of time-integrations of this signal. These
results, obtained through the use of differential algebra and operational
calculus \cite{Mboup-2007journal}, allow to obtain an estimate of the
time-derivative of a particular order in an arbitrary small amount of time
\cite{Fliess-2003kalman}.

Questions arise on how to relate the above to more classical results of
automatic control, and in particular to linear system theory. The present
paper contributes to this discussion by showing that the algebraic
time-derivative estimation method, as presented in \cite{Mboup-2007} and
references therein, can be seen, essentially, as a special case of previously
known state-space results exhibiting a deadbeat property.

After this introduction, we briefly recall in Section \ref{section-algebraic}
the main results of the algebraic time-derivative estimation method. Then,
in Section \ref{section-gramian}, we recall a few results of linear
observability theory and show how in particular the reconstructibility
Gramian can be related to the algebraic method. We end this paper with a few
additional remarks on how to relate further extensions of the algebraic
approach with different areas of control systems theory.

Parts of this study were presented, albeit in German, in Reger and
Jouffroy \cite{RJ-2008German}.

\section{Algebraic time-derivative estimation\label{section-algebraic}}

The algebraic derivative estimation techniques have been presented in
various styles and frameworks, mostly based on abstract algebra and
operational calculus. Because of its practical interest, we recall here only
the main result for a moving-horizon version of the approach (see \cite%
{Mboup-2007} and \cite{Zehetner-2007cca}). However, note that the results
shown in the present paper would also be very easily applicable to earlier
expanding-horizon versions that can be found in \cite%
{Fliess-2005GRETSI-signal} or \cite{Fliess-2004NOLCOS}.

Consider a real-valued, $N$-th degree polynomial function of time
\begin{equation}
y(t)=\sum_{i=0}^{N}\frac{a_{i}}{i!}\,t^{i}  \label{Taylor}
\end{equation}%
where the terms $a_{i}$ are unknown constant coefficients. The goal is to
obtain estimates of the time-derivatives of $y(t)$, up to order $N$.

In \cite{Fliess-2005GRETSI-signal,Fliess-2005GRETSI-image,Neves-2004},
Fliess and co-workers proposed to do so by, roughly speaking, resorting to
algebraic combinations of moving-horizon time-integrations of the available
signal $y(t)$. Let us briefly recall these results in the following theorem
\cite{Mboup-2007,Zehetner-2007cca}.

\begin{theorem}
\label{theorem-algebraic}For all $t\geq T$, the $j$-th order time-derivative
estimate $\hat{y}^{(j)}(t)$, $j=0,1,2,\ldots ,N$, of the polynomial signal $%
y(t)$ as defined in (\ref{Taylor}) satisfies the convolution
\begin{equation}
\hat{y}^{(j)}(t)=\int_{0}^{T}H_j(T,\tau )\,y(t-\tau )\,d\tau \,,\quad
j=0,1,\ldots ,N  \label{a_j_at_t}
\end{equation}%
where the convolution kernel
\begin{align}
& \textstyle H_j(T,\tau )=
\frac{(N + j + 1)!\,(N + 1)!}{T^{N+j+1}}\times \notag\\
& \textstyle\sum\limits_{\kappa_1=0}^{N-j}
\sum\limits_{\kappa_2=0}^{j} \frac{(T-\tau)^{\kappa_1+\kappa_2}\,(-\tau
)^{N-\kappa_1-\kappa_2}}{\kappa_1!\kappa_2!(N-j-\kappa_1)!(j-%
\kappa_2)! (N-\kappa_1-\kappa_2)!(\kappa_1+\kappa_2)!
(N-\kappa_1+1)}  \label{Pi}
\end{align}
depends on the order $j$ of the time derivative to be estimated and on an
arbitrary constant time window length $T>0$.\hfill $\square$
\end{theorem}

\bigskip

For the interested reader, as well as for the sake of completeness, a way to
derive the results of Theorem \ref{theorem-algebraic} is given in Appendix A.

Thus, considering for example the degree-one polynomial%
\begin{equation}
y(t)=a_{0}+a_{1}\,t  \label{pol}
\end{equation}%
applying Theorem \ref{theorem-algebraic} would simply give us the following
first-order time-derivative estimate
\begin{equation}
\hat{\dot{y}}(t)=\int_{0}^{T}\frac{6}{T^{3}}\big(T-2\tau \big)\,y(t-\tau
)\,d\tau\,.  \label{y-dot-est}
\end{equation}%
The effect of the time-integration is obviously to dampen the impact of the
measurement noise on the estimate. Note that this feature can also be used
to filter out noise from the signal $y(t)$ itself, as the zero-order
time-derivative estimator would be
\begin{equation}
\hat{y}(t)=\int_{0}^{T}\frac{2}{T^{2}}\big(2T-3\tau \big)\,y(t-\tau )\,d\tau
\label{y-est}
\end{equation}
as obtained, once again, from Theorem~\ref{theorem-algebraic}.

\section{From deadbeat reconstruction of the state to the algebraic method
\label{section-gramian}}

As will be seen, the above may be related in several ways to more traditional
results of classical linear control theory. To this end, consider now the
following linear time-varying system%
\begin{align}
\dot{\mathbf{x}}(t)& =\mathbf{A}(t)\mathbf{x}(t)  \label{sys-dyn} \\
y(t)& =\mathbf{C}(t)\mathbf{x}(t)  \label{sys-out}
\end{align}%
where $\mathbf{x}(t)\in
\mathbb{R}
^{N+1}$ and $y(t)\in
\mathbb{R}
$. Note that while the form of system (\ref{sys-dyn})-(\ref{sys-out}) was
chosen for the sake of simplicity and ease of presentation, the discussion
of the present section is extendible to systems with multiple inputs and
outputs.

Then let us briefly recall a few elements pertaining to the notion of state
\emph{reconstructibility} \cite{Kalman-1969topics,Brockett-1970book,
Oreilly-1983book}. As noted in Willems and Mitter \cite{Willems-1971tac},
this property has been quite overlooked in the control literature, possibly
because of its equivalence with observability for linear continuous-time
systems. Loosely speaking, we say that system (\ref{sys-dyn})-(\ref{sys-out}%
) is \emph{reconstructible} on $[t_{0},t_{1}]$ if $\mathbf{x}(t_{1})$ can be
obtained from the measurements $y(t)$ for $t\in \lbrack t_{0},t_{1}]$.

A standard way of determining $\mathbf{x}(t_{1})$ can be obtained by first
writing the following expression for the output%
\begin{equation}
y(\tau )=\mathbf{C}(\tau )\,\mathbf{\Phi }(\tau ,t_{1})\,\mathbf{x}(t_{1})
\label{output-eq}
\end{equation}%
where $\mathbf{\Phi }(\tau ,t)$ is the transition matrix of (\ref{sys-dyn}).
Then, left-multiply and integrate (\ref{output-eq}) to get%
\begin{multline}
\int_{t_{0}}^{t_{1}}\mathbf{\Phi }^{\mathrm{T}}(\tau ,t_{1})\,\mathbf{C}^{%
\mathrm{T}}(\tau )\,y(\tau )\,d\tau = \\
\left(\int_{t_{0}}^{t_{1}}\mathbf{\Phi }^{\mathrm{T}}(\tau ,t_{1})\, \mathbf{%
C}^{\mathrm{T}}(\tau )\,\mathbf{C}(\tau )\,\mathbf{\Phi } (\tau
,t_{1})\,d\tau\right)\mathbf{x}(t_{1})  \label{int-eq}
\end{multline}%
Since in eq. (\ref{int-eq}) $\mathbf{x}(t_{1})$ is a constant term with
respect to the integral, it can be isolated, and we finally get, for an
estimate $\hat{\mathbf{x}}(t_{1})$ of $\mathbf{x}(t_{1})$,
\begin{equation}
\hat{\mathbf{x}}(t_{1}):=\mathbf{W}_{\mathrm{r}}^{-1}(t_{0},t_{1})%
\int_{t_{0}}^{t_{1}}\mathbf{\Phi }^{\mathrm{T}}(\tau ,t_{1})\,\mathbf{C}^{%
\mathrm{T}}(\tau )\,y(\tau )\,d\tau  \label{est-gramian}
\end{equation}%
where%
\begin{equation}
\mathbf{W}_{\mathrm{r}}(t_{0},t_{1})=\int_{t_{0}}^{t_{1}}\mathbf{\Phi }^{%
\mathrm{T}}(\tau ,t_{1})\,\mathbf{C}^{\mathrm{T}}(\tau )\,\mathbf{C}(\tau )\,%
\mathbf{\Phi }(\tau ,t_{1})\,d\tau  \label{Gramian}
\end{equation}%
is the \emph{reconstructibility Gramian}.

In treatments of observability in textbooks, developments such as the above
are mostly used, through the observability counterpart of (\ref{Gramian}),
to check whether a system is observable (resp. reconstructible) or not.
However, as noted in \cite[p. 158]{Chen-1999} for the observability case,
expression (\ref{Gramian}) can also be used to actually compute $\mathbf{%
\hat{x}}(t_{1})$ as integration will smooth out high-frequency noise.

The above results are well-known, even if not as much used for state
estimation as linear asymptotic observers are. But the former has the
interesting property of allowing to give an estimate of $\mathbf{x}(t_{1})$
in a \emph{finite} time, whose value is decided by the invertibility of (\ref%
{Gramian}).

\bigskip Interestingly, these two features of the above Gramian-based
estimation -- deadbeat property and time-integration, coincide with those of
algebraic time-derivative estimation.

Let us push the comparison a little further in a simple way by first
noticing that the degree-one polynomial (\ref{pol}) of our example can be put into
state-space phase-variable form with matrices%
\begin{equation}
\mathbf{A}=%
\begin{pmatrix}
0 & 1 \\
0 & 0%
\end{pmatrix}%
, \quad \mathbf{C}=
\begin{pmatrix}
1 & 0%
\end{pmatrix}%
,  \label{pol-matrices}
\end{equation}%
with state $\mathbf{x}(t)=(y(t),\dot{y}(t))^{\mathrm{T}}$ and initial conditions $%
\mathbf{x}(0)=(a_{0},a_{1})^{\mathrm{T}}$.

Then, compute an estimate of $\mathbf{x}(t)$ using (\ref{est-gramian}) and (%
\ref{Gramian}). To do so, use the fact that the matrices in (\ref%
{pol-matrices}) are time-invariant and that $\mathbf{A}^{2}=\mathbf{0}$ to
obtain%
\begin{equation}
\mathbf{\Phi }(\tau ,t_{1})=e^{\mathbf{A}(\tau -t_{1})}=\mathbf{I}+(\tau
-t_{1})\,\mathbf{A}=
\begin{pmatrix}
1 & \tau -t_{1} \\
0 & 1%
\end{pmatrix}%
\end{equation}%
which implies that%
\begin{equation}
\mathbf{C\,\Phi }(\tau ,t_{1})=
\begin{pmatrix}
1 & \tau -t_{1}%
\end{pmatrix}\,.%
\end{equation}%
Letting $t_0=t-T$ (with $T>0$ fixed) and $t_{1}=t$, we then obtain from (\ref%
{Gramian}) the following Gramian%
\begin{equation}
\mathbf{W}_{\mathrm{r}}(t-T,t)=
\begin{pmatrix}
T & -\frac{T^{2}}{2} \\[0.5ex]
-\frac{T^{2}}{2} & \frac{T^{3}}{3}%
\end{pmatrix}%
\end{equation}%
which in turn is used, in combination with (\ref{est-gramian}), to get
\begin{equation}
\hat{\mathbf{x}}(t)=
\begin{pmatrix}
\hat{y}(t) \\
\hat{\dot{y}}(t)%
\end{pmatrix}
=
\begin{pmatrix}
\frac{4}{T} & \frac{6}{T^{2}} \\[0.5ex]
\frac{6}{T^{2}} & \frac{12}{T^{3}}%
\end{pmatrix}
\int_{t-T}^{t}
\begin{pmatrix}
1 \\
\tau -t%
\end{pmatrix}
y(\tau )\,d\tau\,.  \label{pol-est}
\end{equation}%
Hence, similarly to the previous section, an estimate of the derivatives of
a degree-one polynomial can be obtained with time-integrations of the
measured signal, albeit this time using tools from classical control theory.

Note, interestingly, that in this particular example, there is more than a
mere similarity. Indeed, after a simple change of variable $\sigma =t-\tau $
in (\ref{pol-est}), we find exactly the same expressions as (\ref{y-dot-est}%
) and (\ref{y-est}).

The above second-order case can be generalized to obtain the $j$-th
time-derivative of any polynomial simply by specializing $\mathbf{A}(t)$ and
$\mathbf{C}(t)$ in (\ref{sys-dyn})-(\ref{sys-out}) to get a state-space
description of polynomial (\ref{Taylor}), which yields, in phase-variable
form the $N+1$ square matrix

\begin{equation}
\mathbf{A}=%
\begin{pmatrix}
0 & 1 & 0 & \cdots & 0 \\[-1ex]
0 & 0 & 1 & \ddots & 0 \\[-1ex]
\vdots & \vdots & \ddots & \ddots & 0 \\
0 & 0 & 0 & \cdots & 1 \\
0 & 0 & 0 & \cdots & 0%
\end{pmatrix}
\label{deadbeat}
\end{equation}%
and the $N+1$ row vector%
\begin{equation}
\mathbf{C}=%
\begin{pmatrix}
1 & 0 & \cdots & 0%
\end{pmatrix}
\label{deadbeat_C}
\end{equation}%
associated to the state vector $\mathbf{x}(t)=(y(t),\dot{y}%
(t),...,y^{(j)}(t),$ $...,$ $y^{(N)}(t))^{\mathrm{T}}$.

After several steps in line with the previous second-order example, we obtain,
similarly to Section \ref{section-algebraic}, an expression of the $j$-th
time-derivative of a polynomial signal (\ref{Taylor}) based on the
reconstructibility Gramian. This is summarized in the following
theorem.

\bigskip

\begin{theorem}\label{theorem-gramian}
For all $t\geq T$, the $j$-th order time-derivative
estimate $\hat{y}^{(j)}(t)$, $j=0,1,2,\ldots ,N$, of the polynomial signal $%
y(t)$ as defined in (\ref{Taylor}) satisfies the convolution
\begin{equation}
\hat{y}^{(j)}(t)=\int_{0}^{T}G_j(T,\sigma )\,y(t-\sigma )\,d\sigma \,,\quad
j=0,1,\ldots ,N  \label{reconest}
\end{equation}%
where the convolution kernel
\begin{equation}
\textstyle
G_j(T,\tau)=\frac{(N+j+1)!}{T^{j+1}j!(N-j)!} \sum\limits_{k=0}^N
\frac{(-1)^{k}(N+k+1)!}{(j+k+1)(N-k)!(k!)^2} \left(\frac{%
\sigma}{T}\right)^k  \label{Ups}
\end{equation}%
depends on the order $j$ of the time derivative to be estimated and
on an arbitrary constant time window length $T>0$.\hfill $\square$
\end{theorem}

\bigskip

\begin{proof}\/ Broadly speaking, the proof is based on obtaining a closed-form
expression corresponding to equations (\ref{est-gramian}) and
(\ref{Gramian}) for the particular case with matrices (\ref{deadbeat}) and
(\ref{deadbeat_C}).

Since this system is LTI, the corresponding transition matrix results from
the matrix exponential of (\ref{deadbeat}), i.e.
\begin{equation}
e^{\mathbf{A}\,t}=%
\begin{pmatrix}
1 & t & t^{2}/2 & t^{3}/6 & \cdots  & t^{N}/N! \\
0 & 1 & t & t^{2}/2 & \cdots  & t^{N-1}/(N-1)! \\
0 & 0 & 1 & t & \cdots  & t^{N-2}/(N-2)! \\[-1ex]
\vdots  & \vdots  & \vdots  & \vdots  & \ddots  & \vdots  \\
0 & 0 & 0 & 0 & \cdots  & t \\
0 & 0 & 0 & 0 & \cdots  & 1%
\end{pmatrix}
\end{equation}%
which is then used to obtain the state-transition matrix
\begin{equation}
\mathbf{\Phi }(\tau ,t_{1})=e^{\mathbf{A}\,(\tau -t_{1})}\,.
\label{transition-matrix2}
\end{equation}%
Consequently, the entries of the $(N+1)\times (N+1)$ reconstructibility
Gramian matrix (\ref{Gramian}) read
\begin{equation}
\left[ W_{\mathrm{r}}\right] _{ij}\!(t_{0},t_{1})=\!\int_{t_{0}}^{t_{1}}%
\!\textstyle\frac{(\tau -t_{1})^{i+j-2}}{(i-1)!(j-1)!}\,d\tau =\frac{%
-(t_{0}-t_{1})^{i+j-1}}{(i-1)!(j-1)!(i+j-1)}\,.\quad
\label{Qij}
\end{equation}%
In view of (\ref{est-gramian}), the inversion of this Gramian is required. Its
entries are provided in closed-form by Lemma \ref{lemma-inverse} in Appendix
B, that is
\begin{multline}
\left[ W_{\mathrm{r}}^{-1}\right] _{ij}\!(t_{0},t_{1}) =\frac{%
(i-1)!\,(j-1)!\,(i+j-1)}{(t_{1}-t_{0})^{i+j-1}}\times \\
{\binom{N+i}{N+1-j}}{\binom{N+j}{N+1-i}}{\binom{i+j-2}{i-1}}^2.
\label{Qij_inv}
\end{multline}
Hence, by using eq. (\ref{Qij_inv}) regarding the particular form of the
transition matrix (\ref{transition-matrix2}), the $(i+1)$-th component of
$\hat{\mathbf{x}}(t)$ follows from eq. (\ref{est-gramian})
\begin{equation}
\hat{x}_{i+1}(t_{1})=\int_{t_{0}}^{t_{1}}\textstyle\sum\limits_{j=0}^{N}\left[
  W_{\mathrm{r}}^{-1}\right] _{i+1,j+1}(t_{0},t_{1})\,
\frac{(\tau -t_{1})^{j}}{j!}\,y(\tau)\,d\tau \,.
\end{equation}%
In other words, the $j$-th time-derivative estimate of $y(t)$ at time $%
t=t_{1}$ can be obtained from the convolution
\begin{equation}
y^{(j)}(t_{1})=\int_{t_{0}}^{t_{1}}\bar{G}_{j}(t_{1},t_{0},\tau )\,y(\tau
)\,d\tau \,,\quad j=0,1,\ldots ,N  \label{gramian-estimator}
\end{equation}%
where
\begin{multline}
\bar{G}_{j}(t_{1},t_{0},\tau )=\frac{(N\!+\!j\!+\!1)!}{(t_{1}\!-%
\!t_{0})^{j+1}j!(N\!-\!j)!}\times  \\
\sum\limits_{k=0}^{N}\frac{(-1)^{k}(N\!+\!k\!+\!1)!}{(j\!+\!k\!+\!1)(N\!-%
\!k)!(k!)^{2}}\left( \frac{t_{1}\!-\!\tau }{t_{1}\!-\!t_{0}}\right)^{k}.
\end{multline}%
A receding-horizon version of equation (\ref{gramian-estimator}) can then be
obtained as follows: Let $t_{0}=t-T$ (with $T>0$ fixed), and $t_{1}=t$.
Proceed then to the change of variable $\sigma =t-\tau $ to obtain (\ref%
{reconest}) and (\ref{Ups}), which completes the proof of the theorem.
\end{proof}

\bigskip

As might be expected from the above discussion and the second-order example,
it is possible to show an equivalence between the algebraic estimator of
Section \ref{section-algebraic} and the one of Theorem \ref{theorem-gramian},
and this for all $N$. We make this statement precise in the following theorem.

\bigskip

\begin{theorem}
\label{PIUPS_nonu} Let $H_j(T,\tau )$ and $G_j(T,\tau)$ be defined as in (%
\ref{Pi}) and (\ref{Ups}), respectively. Then for $T>0$, $\tau\in[0,T]$ and
$N\in \{0,1,2,\ldots \},$
\begin{equation}
H_j(T,\tau )=G_j(T,\tau ),\quad j\in \{0,1,2,\ldots ,N\}\,.
\end{equation}
\hfill$\square$
\end{theorem}

\bigskip

\begin{proof}\/ Theorem \ref{PIUPS_nonu} follows from Riesz' representation
theorem \cite{Riesz-1955}, which states that for every continuous linear
functional $f$ on a Hilbert space $\mathcal{H}$, a unique $p\in\mathcal{H}$
exists such that
\begin{equation}
f(q)=\langle p,q\rangle \qquad \forall q\in \mathcal{H}\,,
\end{equation}
where $\langle .\,,.\rangle$ denotes the inner product on $\mathcal{H}$.

In order to prepare the ground for applying this theorem, first note that
for parameter $T>0$ fixed, the expressions $H_j(T,\tau)$ and $G_j(T,\tau)$,
given by (\ref{Pi}) and (\ref{Ups}), are polynomials in $\tau$ of degree $N$%
. For $t$ fixed, furthermore $y(t-\tau)$ is a polynomial in $\tau$ of degree
$N$ which in view of (\ref{Taylor}) consequently spans $\mathcal{H}_N$, i.e.
the Hilbert space of degree $N$ polynomials equipped with the real-valued
inner product
\begin{equation}
\langle p,q\rangle :=\int_{0}^{T}p(\tau )q(\tau )\,d\tau \,,\quad p,q\in%
\mathcal{H}_N\,.
\end{equation}
Hence, for $T>0$ fixed, $H_j(T,\tau)\in\mathcal{H}_N$ and $G_j(T,\tau)\in%
\mathcal{H}_N$. Moreover, letting $q(\tau):=y(t-\tau)$ with fixed $t\ge T$
we have that $q\in\mathcal{H}_N$.

In accordance with (\ref{a_j_at_t}) and (\ref{reconest}), let
\begin{equation}
f_{H_j}(q):=\int_{0}^{T}H_{j}(T,\tau )\,q(\tau )\,d\tau
\end{equation}%
and%
\begin{equation}
f_{G_j}(q):=\int_{0}^{T}G_{j}(T,\tau )\,q(\tau )\,d\tau
\end{equation}%
for $j=0,1,2,\ldots ,N$.

Thus, Theorems \ref{theorem-algebraic} and \ref{theorem-gramian}
imply that for any $q\in\mathcal{H}_N$
\begin{equation}
f_{H_j}(q)=f_{G_j}(q)\,,\quad j=0,1,2,\ldots,N\,.
\end{equation}%
Since $H_j(T,\tau)\in\mathcal{H}_N$ and $G_j(T,\tau)\in\mathcal{H}_N$, for $%
T>0$ fixed, the uniqueness of $p$ in Riesz' theorem shows that
\begin{equation}
H_j(T,\tau )\equiv G_j(T,\tau )
\end{equation}
for $j=0,1,2,\ldots ,N$, under the assumptions of Theorem \ref{PIUPS_nonu}.
\end{proof}
Note that other proofs of the previous theorem are also possible. For example,
a somewhat more component-wise proof, based on modern computer algebra proof
techniques \cite{WZ92a}, is presented in \cite{RJ-2007Report} by showing
specifically how the terms in (\ref{Pi}) relate to those of (\ref{Ups}).

\section{Additional remarks}

In addition to the main result of Section \ref{section-algebraic}, Fliess et
al. proposed several extensions or modifications, several of which having
also connections with different areas of control systems. Let us briefly
consider some of them in the few following remarks.

For instance, note that an \emph{expanding-horizon} version of the algebraic
method was first introduced in \cite{Fliess-2003kalman}, which would
correspond to let $t_0=0 $ and $t_1=t$ in the reconstructibility Gramian
perspective. In this case, an equivalence similar to Theorem \ref{PIUPS_nonu}
can still be obtained. Furthermore, note that, interestingly, letting $%
\mathbf{S}(t):=\mathbf{W}_{\mathrm{r}}(0,t)$, and differentiating
respectively $\mathbf{S}(t)$ and the product $\mathbf{S}(t)\,\hat{\mathbf{x}}%
(t)$ with respect to time using a few standard manipulations, we obtain
\begin{equation}
\dot{\mathbf{S}}(t)=-\mathbf{A}^{\mathrm{T}}(t)\mathbf{S}(t)-\mathbf{S}(t)%
\mathbf{A}(t)+\mathbf{C}^{\mathrm{T}}(t)\mathbf{C}(t)\vspace{-1ex}
\end{equation}
and
\begin{equation}
\dot{\hat{\mathbf{x}}}(t)=\left( \mathbf{A}(t)-\mathbf{S}^{-1}(t)\mathbf{C}^{%
\mathrm{T}}(t)\mathbf{C}(t)\right) \hat{\mathbf{x}}(t)+\mathbf{S}^{-1}(t)%
\mathbf{C}^{\mathrm{T}}(t)y(t)
\end{equation}
which draw strong similarities with the information form of the continuous-time
Kalman filter \cite{Kwon-1999,Kaminski-1971tac} for system (\ref{sys-dyn})-(%
\ref{sys-out}) with additive noise $v(t)\in
\mathbb{R}
$ of identity covariance, $\mathbf{R}=\mathbf{I}$, on the measurement
equation (\ref{sys-out}). This in turn shows that, thanks to a simple
modification of Theorem \ref{PIUPS_nonu} for expanding horizons, links with
optimal estimation could be obtained even though the derivations and
motivations for the algebraic method are clearly different (see in
particular \cite{Fliess-2003kalman}).

As another example, one could consider the state estimation problem for a MIMO
time-invariant system with inputs:
\begin{align}
\dot{\mathbf{x}}(t)& =\mathbf{A}\mathbf{x}(t)+\mathbf{B}\mathbf{u}(t)
\label{sys-dyn-byrski} \\
\mathbf{y}(t)& =\mathbf{C}\mathbf{x}(t)  \label{sys-out-byrski}
\end{align}%
For such systems, Byrski et al. derived a so-called {\it moving window
  observer}
\cite{Byrski-1984cac},\cite{Fuksa-1984tac},\cite{Byrski-1993ecc}. In order to
briefly sketch the result, let
\begin{equation}
\mathbf{\Psi}(t)=e^{\mathbf{\Omega}t}
\end{equation}
 with
\begin{equation}
\mathbf{\Omega}=
\begin{pmatrix}
\mathbf{A} & \mathbf{B}\mathbf{B}^{\mathrm{T}}\\
\mathbf{C}^{\mathrm{T}}\mathbf{C} & -\mathbf{A}^{\mathrm{T}}
\end{pmatrix}\,.
\end{equation}
Assume that the output and the input of system (\ref{sys-dyn-byrski}) and
(\ref{sys-out-byrski}) may be used for the reconstruction of the state, thus
consequently, we may allow for input- and output-sided deterministic
disturbances, bounded in an $\mathcal{L}_2$-norm sense. Moreover, assume that
the pair $(\mathbf{C},\mathbf{A})$ is observable. Then the {\it moving window
  observer} that minimizes the estimation error $\hat{\mathbf{x}}-\mathbf{x}$
on the moving fixed time horizon $[t-T,t]$ is given by
\begin{equation}
\hat{\mathbf{x}}(t)\!=\!\!\int_0^T\!\!\!\!\mathbf{G}_1
(T,T\!-\tau)\mathbf{y}(t-\tau)\,d\tau+
\!\int_0^T\!\!\!\!\mathbf{G}_2(T,T\!-\tau)\mathbf{u}(t-\tau)\,d\tau
\label{est-byrski}
\end{equation}
where
\begin{align}
\mathbf{G}_1(T,t)= & e^{\mathbf{A}T}\left(
\int_0^T\mathbf{\Psi}_{11}^{\mathrm{T}}(\tau)\,\mathbf{C}^{\mathrm{T}}
\mathbf{C}\,e^{\mathbf{A}\tau}d\tau
\right)^{-1}
\mathbf{\Psi}_{11}^{\mathrm{T}}(t)\,\mathbf{C}^{\mathrm{T}}\label{G1}\\
\mathbf{G}_2(T,t)= & e^{\mathbf{A}T}\left(
\int_0^T\mathbf{\Psi}_{11}^{\mathrm{T}}(\tau)\,\mathbf{C}^{\mathrm{T}}
\mathbf{C}\,e^{\mathbf{A}\tau}d\tau
\right)^{-1}
\mathbf{\Psi}_{21}^{\mathrm{T}}(t)\,\mathbf{B}\label{G2}
\end{align}
In the case of an input-free SISO system of the particular form
(\ref{deadbeat}) and (\ref{deadbeat_C}), we have that
$\mathbf{\Psi}_{11}(t)=e^{\mathbf{A}t}$ and matrix $\mathbf{G}_2$ vanishes. As
a consequence, algebraic derivative estimation may be seen as a very special
particularization of equation (\ref{est-byrski}).


In an other extension presented in \cite{Neves-2004}, the authors propose to
further reduce the impact of measurement noise on the estimates by using
additional integrations. This is also possible with the Gramian point-of-view
as both sides of (\ref{int-eq}) can easily be time-integrated several
additional times with respect to $t_{0}$, as opposed to only once to obtain
$\mathbf{x}(t_{1})$ -- in fact, even filter operations with respect to the
variable $t_{0}$ can be applied on both sides of (\ref{int-eq}), so as to
generate a variety of further estimators. Once again, an equivalence between
this result of the algebraic approach and a particularization of a
reconstructibility perspective can be obtained. More generally, we can for
example insert in (\ref{int-eq}) another kernel $\lambda (\tau ,t_{0})$ as
follows
\begin{equation}
\hat{\mathbf{x}}(t_{1}):=\mathbf{W}_{\mathbf{\lambda }}^{-1}(t_{0},t_{1})%
\int_{t_{0}}^{t_{1}}\lambda (\tau ,t_{0})\mathbf{\Phi }^{\mathrm{T}}(\tau
,t_{1})\,\mathbf{C}^{\mathrm{T}}(\tau )\,y(\tau )\,d\tau
\end{equation}
where
\begin{equation}
\mathbf{W}_{\mathbf{\lambda }}(t_0,t_1)\!=\!\!\int_{t_0}^{t_1}\!\!\lambda
(\tau ,t_{0})\mathbf{\Phi }^{\mathrm{T}}(\tau ,t_{1})\,\mathbf{C}^{\mathrm{T}%
}(\tau )\,\mathbf{C}(\tau )\,\mathbf{\Phi }(\tau ,t_{1})\,d\tau ,
\end{equation}%
this to obtain the desired response with respect to measurement noise.

Finally, and although it is clearly beyond the scope of the present paper,
note that because of the convolution form of algebraic estimation (\ref%
{a_j_at_t}), the latter can also be connected with Finite-Impulse Response
(FIR) differentiators, on which numerous studies and results were published
(see \cite{Khan-1999iee}, \cite{Tseng-2002} and references therein), with
the minor difference that these differentiators are usually described in a
discrete-time framework, although it is clear that a comparison similar to
the present paper could also be carried out in discrete-time.

In particular, it might be of interest to compare the latest extension of
the algebraic estimation approach, where time-delays are considered to
improve the results, together with FIR differentiator designs considering
the same issue that have been proposed over the past few years (see for
example \cite{Vainio-1999tim}\ and \cite{Samadi-2007}).

\section*{Acknowledgments}

The authors would like to express their gratitude to H\aa kan Hjalmarsson
whose comments and suggestions were very helpful in improving the present
paper. Johann Reger thanks Peter Caines and Jessy Grizzle for valuable
discussions. The work of Johann Reger was partially supported by fellowships
of the German Academic Exchange Service (DAAD), grant D/07/40582, and by Max
Planck Institute for Dynamics of Complex Technical Systems in Magdeburg,
Germany.

\section*{Appendix}

\subsection{Proof of Theorem \protect\ref{theorem-algebraic}}

The following proof resorts to standard techniques from operational calculus.
To this end, we rephrase eq. (\ref{Taylor}) in the Laplace domain as
\begin{equation}
Y(s)=\sum_{i=0}^{N}\frac{y^{(i)}(0)}{s^{i+1}}\,,  \label{YN}
\end{equation}%
where the coefficients $a_{i}$ are identified with $y^{(i)}(0)$. In order to
single out a particular term, $y^{(j)}(0)$, first multiply (\ref{YN}) by $%
s^{N+1}$,
\begin{equation}
s^{N+1}\,Y(s)=\sum_{i=0}^{N}y^{(i)}(0)\,s^{N-i}\,,  \label{operation1}
\end{equation}%
which results in a polynomial form in $s$ on the right side of (\ref%
{operation1}). To eliminate the terms $y^{(j+1)}(0),\ldots ,y^{(N)}(0)$,
differentiate (\ref{operation1}) $N-j$ times with respect to $s$ (see \cite%
{Fliess-2003mex} for a first presentation of the idea). This yields
\begin{equation}
\frac{d^{N-j}}{ds^{N-j}}\left( s^{N+1}Y(s)\right) =\sum_{i=0}^{j}y^{(i)}(0)\,%
\frac{(N-i)!}{(j-i)!}\,s^{j-i}\,.  \label{operation2}
\end{equation}%
In the next step, we proceed to a similar treatment to eliminate the
remaining constant terms $y^{(0)}(0)$, $y^{(1)}(0)$, $\ldots $, $y^{(j-1)}(0)
$. But before doing so, premultiply (\ref{operation2}) by $1/s$, i.e.
\begin{multline}
\frac{1}{s}\frac{d^{N-j}}{ds^{N-j}}\left(s^{N+1}Y(s)\right)=\\
=\frac{(N-j)!}{s}y^{(j)}(0)+\sum_{i=0}^{j-1}y^{(i)}(0)\frac{(N-i)!}{%
(j-i)!}s^{j-i-1}
\end{multline}%
which is done to prevent $y^{(j)}(0)$ from cancelation due to a $j$-fold
differentiation with respect to $s$. Indeed, the latter operation finally
gives
\begin{equation}
\frac{d^{j}}{ds^{j}}\!\!\left( \frac{1}{s}\frac{d^{N-j}}{ds^{N-j}}\!\left(
s^{N+1}Y(s)\right) \!\!\right) \!\!=\!\frac{(-1)^{j}\,j!\,(N\!-\!j)!}{s^{j+1}%
}y^{(j)}(0).  \label{isolated}
\end{equation}
This equation could readily be transformed back into the time domain.
However, the left side of (\ref{isolated}) contains the monomial $s^{N}$,
i.e. an $N$-fold differentiation with respect to time in the time domain,
meaning if a high-frequency noise is corrupting $y(t)$, the former would be
amplified as a result. Note that a similar idea can also be found in \cite[%
p.17--18]{Sontag-1998book}. In order to avoid the explicit use of these time
derivatives, premultiply (\ref{isolated}) with $1/s^{N+1}$, thus implying
that $y(t)$ will be integrated at least one time. Therefore, we obtain
\begin{equation}
\frac{1}{s^{N+1}}\frac{d^{j}}{ds^{j}}\!\!\left(\!\frac{1}{s}\frac{d^{N-j}}{%
ds^{N-j}}\!\left(s^{N+1}Y(s)\!\right)\!\!\right) \!\!=\!\frac{%
(-1)^{j}j!\,(N\!-\!j)!}{s^{N+j+2}}y^{(j)}(0)  \label{vorleibniz}
\end{equation}%
where it can been seen that the term $y^{(j)}(0)$ depends only on a finite
number of operations on the signal $Y(s)$, as shown in \cite%
{Mboup-2007,Zehetner-2007cca}.

Before performing the backward transform into the time-domain, rearrange the
left side terms of (\ref{vorleibniz}) using Leibniz' formula for the
differentiation of products twice. This results in
\begin{multline}
\!\!\!\!\!\!\frac{1}{s^{N+1}}\frac{d^{j}}{ds^{j}}\!\left(\!\frac{1}{s}
\frac{d^{N-j}}{ds^{N-j}}\!\left(s^{N+1}Y(s)\!\right) \!\!\right)
\!\!=\!\!\sum\limits_{\kappa _{1}=0}^{N-j}\sum\limits_{\kappa _{2}=0}^{j}\!\!{%
\binom{N\!-\!j}{\kappa _{1}}}\!\!{\binom{j}{\kappa _{2}}\times } \\
\frac{(N\!+\!1)!}{(N\!-\!\kappa _{1}\!-\!\kappa _{2})!\,(N\!-\!\kappa
_{1}\!+\!1)}\frac{1}{s^{\kappa _{1}+\kappa _{2}+1}}\frac{d^{N-\kappa
_{1}-\kappa _{2}}}{ds^{N-\kappa _{1}-\kappa _{2}}}\,Y(s)
\end{multline}%
which, in view of the right hand side of (\ref{vorleibniz}), implies in turn
\begin{multline}
\frac{1}{s^{N+j+2}}\,y^{(j)}(0)=\frac{(-1)^{j}}{j!\,(N\!-\!j)!}%
\sum\limits_{\kappa _{1}=0}^{N-j}\sum\limits_{\kappa _{2}=0}^{j}\!{\binom{%
N\!-\!j}{\kappa _{1}}}\!\!{\binom{j}{\kappa _{2}}}\times  \\
\frac{(N\!+\!1)!}{(N\!-\!\kappa _{1}\!-\!\kappa _{2})!\,(N\!-\!\kappa
_{1}\!+\!1)}\frac{1}{s^{\kappa _{1}+\kappa _{2}+1}}\frac{d^{N-\kappa
_{1}-\kappa _{2}}}{ds^{N-\kappa _{1}-\kappa _{2}}}\,Y(s)\,.
\label{isolated-LT}
\end{multline}%
Eq. (\ref{isolated-LT}) is now transformed back into the time domain. Using
the following inverse Laplace transform formulae%
\begin{equation}
\mathcal{L} ^{-1}\ \left[ \frac{1}{s^{i+1}}\frac{d^{j}}{ds^{j}}Y(s)\right] =\
\int_{0}^{t}\frac{(t-\tau )^{i}(-\tau )^{j}}{i!}\,y(\tau )\,d\tau
\end{equation}

\noindent we obtain
\begin{equation}
\hat{y}^{(j)}(0)=\int_{0}^{t}H_{j}(t,\tau )\,y(\tau )\,d\tau \,,\quad
j=0,1,\ldots ,N
\end{equation}%
with
\begin{align}
  & \textstyle H_{j}(t,\tau)=\frac{(N+j+1)!\,(N+1)!\,(-1)^{j}}{t^{N+j+1}}%
  \times\notag\\
  & \!\textstyle\sum\limits_{\kappa
    _{1}=0}^{N-j}\sum\limits_{\kappa _{2}=0}^{j}\!\frac{(t-\tau )^{\kappa
      _{1}+\kappa _{2}}\,(-\tau )^{N-\kappa _{1}-\kappa _{2}}}{\kappa
    _{1}!\kappa _{2}!(N-j-\kappa _{1})!(j-\kappa
    _{2})!(N-\kappa _{1}-\kappa _{2})!(\kappa _{1}+\kappa
    _{2})!(N-\kappa _{1}+1)}
\end{align}
The results obtained above thus give an estimate $\hat{y}^{(j)}(t)$ at
time $t=0$ from the polynomial signal $y$, see (\ref{Taylor}), taken on the
interval $[0,t]$. In order to get a moving-horizon and causal version of
these results, first replace $t$ with $-T$, where $T$ is a positive constant
\cite{Fliess-2005GRETSI-signal,Fliess-2005GRETSI-image} and simplify using
the fact that
\begin{equation}
(-1)\ H_{j}(-T,-\tau )=(-1)^{j}\ H_{j}(T,\tau )\,.
\end{equation}%
Finally, by shifting the $y$-values by $t$, Theorem \ref{theorem-algebraic}
is immediate.\hfill $\blacksquare $

\subsection{Lemma for the Proof of Theorem \ref{theorem-gramian}}

\begin{lemma}[Inverse of $\mathbf{W}_{\mathrm{r}}(t_{0},t_{1})$]
\label{lemma-inverse}Let the entries of the matrix $\mathbf{W}_{\mathrm{r}%
}(t_{0},t_{1})$ be given as in (\ref{Qij}). The entries of its inverse
are
\begin{multline}
\left[ W_{\mathrm{r}}^{-1}\right] _{ij}(t_{0},t_{1}) =
\frac{(i-1)!\,(j-1)!\,(i+j-1)}{(t_{1}-t_{0})^{i+j-1}}\times\\
{\binom{N+i}{N+1-j}}{\binom{N+j}{N+1-i}}{\binom{i+j-2}{i-1}}^2.
\end{multline}%
\hfill $\square$
\end{lemma}

\bigskip

\begin{proof}\/ In light of equation (\ref{Qij}), first, left- and
  right-multiply $\mathbf{W}_{\mathrm{r}}(t_{0},t_{1})$ with a diagonal matrix
  $\mathbf{M}$ whose entries are
\begin{equation}
M_{ij}=\frac{\,(i-1)!\,}{(t_{0}-t_{1})^{i}}\delta _{ij}
\end{equation}%
where $\delta _{ij}$ is the Kronecker delta. Then, proceed with computing
the following matrix product in component form as
\begin{align}
& \lbrack (t_1\!-\!t_0)\,\mathbf{M}\,\mathbf{W}_{\mathrm{r}}(t_{0},t_{1})\,%
\mathbf{M}]_{ij}  \notag \\
& \textstyle =(t_1-t_0)\sum\limits_{k,l=1}^{N+1} M_{ik}
\left[ W_{\mathrm{r}}\right]_{kl}(t_{0},t_{1})M_{lj}\notag \\
& \textstyle=(t_1-t_0)\sum\limits_{k,l=1}^{N+1}
\frac{(i-1)!}{(t_0-t_1)^i} \delta_{ik} \frac{-(t_0-t_1)^{k+l-1}}{%
(k-1)!(l-1)!(k+l-1)} \frac{(l-1)!}{(t_0-t_1)^l}%
\delta_{lj}  \notag \\
&  =\frac{1}{i+j-1}
\end{align}%
whose result can be recognized as the entries of an $(N+1)\times$ $(N+1)$
Hilbert matrix, hereafter denoted $\mathbf{H}$. The entries of the inverse
of $\mathbf{H}$ are known to be \cite{Savage-1954}
\begin{equation}
\left[ H^{-1}\right]_{ij}=(-1)^{i+j}\,(i+j-1)
\textstyle
{\binom{N+i}{N+1-j}}
{\binom{N+j}{N+1-i}}
{\binom{i+j-2}{i-1}}^2\quad
\end{equation}%
and by computing
\begin{equation}
\mathbf{W}_{\mathrm{r}}^{-1}(t_0,t_1)=(t_1\!-\!t_0)\,\mathbf{M}\,\mathbf{H}%
^{-1}\,\mathbf{M}
\end{equation}%
we obtain (\ref{Qij_inv}), which completes the proof of the Lemma.
\end{proof}

\vspace{8.5mm}

\end{document}